\documentclass[preprint]{elsarticle}
\usepackage{amsmath,amsthm, amssymb, amsfonts}
\DeclareMathOperator{\Lie}{\mathrm{Lie}}
\DeclareMathOperator{\cA}{\mathcal{A}}
\DeclareMathOperator{\cO}{\mathcal{O}}
\DeclareMathOperator{\dom}{\partial_\omega}
\DeclareMathOperator{\dth}{\partial_\theta}
\DeclareMathOperator{\Ve}{Vec}
\DeclareMathOperator{\z}{\mathbf{0}}

\newcommand{\bbox}{\hfill $\Box$}

\newcommand{\F}{\mathcal{F}}
\newcommand{\N}{\mathbb{N}}
\newcommand{\p}{\partial}
\newcommand{\R}{\mathbb{R}}
\newtheorem{corollary}{Corollary}
\newtheorem{lemma}{Lemma}
\newtheorem{proposition}{Proposition}
\newtheorem{theorem}{Theorem}
\newproof{pf}{Proof}
\begin{document}
%%%%%%%%%%%%%%%%%%%%%%%%%%%%%%%%%%%%%%
\begin{frontmatter}
%%%%%%%%%%%%%%%%%%%%%%%%%%%%%%%%%%%%%%
\title{Controllability of a planar inverted pendulum on a cart}
\author{M.Isabel Caiado\fnref{fn1}}
\ead{icaiado@math.uminho.pt}
\address{Centro de Matem\'{a}tica, Universidade do Minho, Campus de Gualtar, 4710-057 Braga, Portugal}
\fntext[fn1]{Work partially supported by Portuguese Foundation for Science and Technology (FCT) under CMAT}
%%%%%%%%%%%%%%%%%%%%%%%%%%%%%%%%%%%%%%
\begin{abstract}
Using tools and techniques from  geometric control theory, namely Hermann-Nagano theorem, we prove controllability result for a mechanical system consisting of a planar double inverted pendulum fixed to a wheeled cart that can move linearly along a horizontal track.
\end{abstract}
%%%%%%%%%%%%%%%%%%%%%%%%%%%%%%%%%%%%%%
\end{frontmatter}
%\maketitle
%%%%%%%%%%%%%%%%%%%%%%%%%%%%%%%%%%%%%%
\section{Introduction}
The inverted pendulum is a typical example of nonlinear control system which is intensively  studied.
The problem of controllability and stabilization  of inverted pendulum on a cart has been addressed by several authors using different strategies.
A series of papers  (\cite{AstFur00, BloLeoMar, Shir01}) has been published  on energy strategies. More recently, Mason {\it et al} in \cite{MaBrPiccoli} obtained qualitative results on the global structure of the time optimal trajectories of the planar pendulum
on a cart.

The present paper aims at a different target: verifying controllability of the (non-linear) mechanical system by application of tools from geometric control theory. 

The paper is organized as follows. After the Introduction, in the second section we describe the mechanical system consisting in two linked planar inverted pendulums on a cart (PIDP) and deduce the equations of the dynamics. The third section recalls some definitions and results from geometrical control theory. In the next section several results are proved in particular the main result which establishes the controllability of PIDP.  The proof of controllability  is mainly based on the existence of a bracket-generating family of vector fields for PIDP which the Lie algebra has constant dimension over each orbit of the family through a point.  The last section contains some remarks on the choice of the bracket-generating family used and on some particular mechanical systems which controllability could not be proved.
%%%%%%%%%%%%%%%%%%%%%%%%%%%%%%%%%%%%%%
\section{Inverted double pendulum on a cart}
%%%%%%%%%%%%%%%%%%%%%%%%%%%%%%%%%%%%%%
We consider a mechanical system which is a double mathematical inverted pendulums
(i.e. a double pendulum in a gravitational field without friction and tension). Each pendulum
is modeled by a mass point (the bob of mass $m_i$) and a massless beam of length $r_i$. The
second pendulum is attached to the bob of the first one. We neglect the axial rotation
of the beams, so there is one degree of freedom for each pendulum. The pivot of the first pendulum is fixed to a wheeled cart that can move linearly along a horizontal track.

We assume that controlled acceleration $u(\cdot)$ is applied to move the cart in horizontally direction in order to balance
the two linked inverted pendulums on the cart. Therefore the motion occurs in the vertical plane.
Consider a coordinate system with the origin at the pivot of the first pendulum. For $i\in \{1,2\}$ let $\theta_i$ be the angle between each pendulum and the positive part of vertical
axis at time $t$, i.e. $\theta_i\equiv\theta_i(t)$ with $\theta_i:[0, \infty[\, \longrightarrow [-\pi,\pi[\,; i\in \{1,2\}$.
On each mass there is an actuating force resulting from the acceleration $u$. Letting
$\theta(\cdot)=(\theta_1(\cdot),\theta_2(\cdot))$, the system's kinetic energy is given by
\begin{equation*}\label{Ener_Cin}
  T(\theta,\dot\theta,t)= \frac{1}{2}\,(m_1+m_2)\,r_1^2\,\dot
  \theta_1^2+m_2\,r_1\,r_2\,\dot\theta_1\,\dot\theta_2\,\cos(\theta_1-\theta_2)+\frac{1}{2}\,m_2\,r_2^2\,\dot
  \theta_2^2
  \end{equation*}
while the potential energy is $U(\theta,t)=g\,[(m_1+m_2)\,r_1\,\cos\theta_1+m_2\,r_2\,\cos\theta_2].$

Using Euler-Lagrange equations, we obtain the dynamical equations.
Since the system is conservative, we get
\begin{equation}
\label{eq: sis Lagrangiano}
\frac{\p L}{\partial \theta_i}-\frac{d}{dt}\frac{\p L}{\p
\dot{\theta}_i}=m_i\,u(t),\qquad i\in \{1,2\},
\end{equation}
where $L=L(\theta,\dot\theta,t)= T(\theta,\dot\theta,t)-U(\theta,t)$
and $u:[0, \infty[\, \longrightarrow \R$ is the control. This is a  system of second order ordinary differential equations in
$\theta_i;i\in \{1,2\}$.
Introducing the notation $r_1\,r_2[m_1+m_2\,\sin^2(\theta_1-\theta_2)]=\Delta(\theta)> 0$,
$\dot \theta_i=\omega_i;\,i\in \{1,2\}$ and $z=(\theta_1,\theta_2, \omega_1,
\omega_2)^T$, the system can be rewritten as a system of first order differential equations
\begin{equation}
\label{eq:sistCont}
  \dot z=f(z)+h(z)u(t),
\end{equation}
where
\begin{align}
f(z)&=\frac{1}{\Delta(\theta)}\begin{pmatrix}
    \omega_1 \,\Delta(\theta)\\
    \omega_2 \,\Delta(\theta)\\
    -g\,r_2\,[m_2\,\sin\,\theta_2\,\cos(\theta_1-\theta_2)-(m_1+m_2)\,\sin\,\theta_1] \\
    -g\,r_1\,(m_1+m_2)\,[\cos(\theta_1-\theta_2)\,\sin\theta_1-\sin\theta_2]
  \end{pmatrix} \label{f}\\
  &+\frac{\sin(\theta_1-\theta_2)}{\Delta(\theta)}\begin{pmatrix}
    0\\
    0\\
    -m_2\,r_2\,[r_1\,\cos(\theta_1-\theta_2)\,\omega_1^2+r_2\,\omega_2^2] \\
    r_1\,[r_1\,(m_1+m_2)\,\omega_1^2+r_2\,m_2\,\cos(\theta_1-\theta_2)\,\omega_2^2] \
  \end{pmatrix} \nonumber\\
  \intertext{and}
  h(z)&=  \frac{m_2}{\Delta(\theta)}
  \begin{pmatrix}
    0 \\
    0 \\
    r_2\,[r_1\,m_2\,\cos(\theta_1-\theta_2)-r_2\,m_1] \\
    r_1\,[r_2\,m_1\cos(\theta_1-\theta_2)-r_1\,(m_1+m_2)]
  \end{pmatrix}.\label{h}
\end{align}

For simplicity, we use the following notation
\begin{align}\label{eq: notacao 1}
\Delta(\theta)\,f(z)=\begin{pmatrix}
   \Omega_1(z)\\
    \Omega_2(z)\\
   a_1(z)\\
    a_2(z)\
  \end{pmatrix}, && \Delta(\theta)\,h(z)=\begin{pmatrix}
    \z_2\\
    b_1(\theta)\\
    b_2(\theta)
  \end{pmatrix}, &&
    \frac{\partial}{\partial\theta}=\partial_{\theta}, &&
    \frac{\partial}{\partial\omega}={\partial_{\omega}},
  \end{align}
  \begin{align}\label{eq: notacao 2}
\Omega(z)=\begin{pmatrix}
   \Omega_1(z)\\
    \Omega_2(z)\\
     \end{pmatrix}, &&
 a(z)=\begin{pmatrix}
   a_1(z)\\
    a_2(z)\
  \end{pmatrix}, &&b(\theta)=\begin{pmatrix}
    b_1(\theta)\\
    b_2(\theta)\end{pmatrix},
  \end{align}
\begin{align}\label{eq: notacao 3}
\bar \Omega (\theta)=-\Delta(\theta)b(\theta), &&\bar a^T (z)=\Omega^T(z)\dth b-b^T(\theta)\dom a,
\end{align}
and
\begin{align}\label{eq: notacao 4}
   \bar b^T(\theta)=b^T(\theta)[2
\Delta(\theta) \dth b- \dom(b^T(\theta)\dom a) ]
\end{align}
where the superscript $T$ stands for transposition and $\z_2$ is the zero two dimensional vector. {We assume that  mechanical system parameters satisfy
\begin{align}\label{eq: nao BG}
    &\frac{m_1}{m_2}\neq\frac{r_1\,(r_1\pm r_2)}{r_2^2-r_1\,(r_1\pm r_2)} && \text{and} &&
    &\frac{m_1}{m_2}\neq\frac{r_2^2-r_1^2}{r_1^2}.
\end{align}
}

To study controllability of system \eqref{eq:sistCont}-\eqref{eq: nao BG}  we will use some techniques from geometric control theory for systems with recurrent drift.  Next section  collects some definitions as well as some classical results.
%%%%%%%%%%%%%%%%%%%%%%%%%%%%%%%%%%%%%%
\section{Classical results from geometric control theory}
%%%%%%%%%%%%%%%%%%%%%%%%%%%%%%%%%%%%%%
Let $V$ be a smooth manifold and $\Ve V$ the set of all smooth vector fields on $V$.
Let $X\in \Ve V$ be a complete vector field. A point $p\in V$ is called Poisson stable for $X$ if for any $t>0$ and any neighbourhood $U$ of $p$ there exists a point $q\in U$ and a time $t'>t$ such that  $q\circ P_{t'}\in\,U$, where $P_t$ is the flow generated by  $X$.
Poisson stability is associated with the vector field defining the system. According to Poisson stability all trajectories cannot leave a neighbourhood of a Poisson stable point forever, some of them must return to this neighbourhood for arbitrarily large times.
Note that, if a trajectory $p\circ P_{t}$ is periodic, then $p$ is Poisson stable for vector field $X$. A complete vector field $X\in \Ve V$ is Poisson stable if all points of $V$ are Poisson stable for $X$.   Next theorem characterizes such sets see  \cite{Bonard, Lobry}.

\begin{theorem}
\label{thm: Poicare}
Let $V$ be a smooth manifold with a volume form $\mathrm{Vol}$. Let a vector field $X \in \Ve \, V$
be complete and its flow $P_t$ preserve volume. Let  $W\subset\,V, W \subset
\overline{\mathrm{int}\, W},$ be a subset of finite volume, invariant for $X$:
\begin{equation*}
    \mathrm{Vol}(W)< \infty, \quad W\circ P_t\subset W,\quad t>0.
\end{equation*}
Then all points of $W$ are Poisson stable for $X$.
\end{theorem}

\medskip

Let
$$\F=\{X_1,\dots,X_k\}\subset \Ve V$$
be a family of complete differentiable vector fields defined on the manifold $V$.
The Lie algebra generated by $\F$ is defined by
\begin{equation*}
  \Lie(\F)=\mathrm{span} \{[X_1,[\dots[X_{k-1},X_k]\dots]]\,:\,k\in  \mathbb{N}, \,X_1,\dots,X_k \in\F,\},
\end{equation*}
where $[\cdot,\cdot]$ stands for Lie bracket.

For any $p\in V$, we denote by $\Lie_p(\F)$ the set of all tangent vectors $p\circ\,X$ to $V$ with $X\in \Lie\,(\F)$
\begin{equation*}
 \Lie_p(\F)=\{p\circ\,X\,:\,X \in \Lie\,(\F)\}\subset
 T_p\,V
\end{equation*}
where  $T_p\,V$ is the tangent space to $V$ at $p$.

The orbit of the family $\F$ through a point $p\in V$ is the set{
\begin{equation*}
  \cO_p(\F)=\{p\circ e^{t_1X_1}\circ \dots\circ e^{t_iX_i}:\,i\in
  \mathbb{N}, t_1, \dots,t_i \in \R, X_1, \dots, X_i\in \F \},
\end{equation*}}
where $e^{t_jX_j}$ is the flow generated by vector field $X_j$, i.e. $\cO_p(\F)$ is the set of all attainable points from $p$ composing flows generated by vector fields in $\F$.
Therefore, in an orbit $\cO_p(\F)$, one is allowed to move along vector fields $X_j$ both forward and backwards (with reverse time). If only forward motion is permitted the set obtained is called attainable set{ \begin{equation*}
  \cA_p(\F)=\{p\circ e^{t_1X_1}\circ \dots \circ e^{t_iX_i}:\,i\in
  \mathbb{N}, t_1, \dots, t_i\geq 0, X_1, \dots, X_i\in \F\}.
\end{equation*}}

For analytic manifolds, the following result, known as Hermann-Nagano theorem, holds (see \cite{Jur} for further information).

\begin{theorem}
\label{HerNag}
Let $V$ be an analytic manifold and $\F$ a family of analytic vector fields on $V$. Then
\begin{enumerate}
  \item[(a)] each orbit of $\F$ is an analytic submanifold of $V$, and
  \item[(b)] if $\cO_p(\F)$ is an orbit of $\F$ through $p\in V$, the tangent space of  $\cO_p(\F)$ at $q$ is given by $\Lie_q(\F)$.
      In particular, the dimension of $\Lie_q(\F)$ is constant as $q$ varies over $\cO_p(\F)$.
\end{enumerate}
\end{theorem}

\medskip

It remains to introduce two definitions in order to state the results which will allow us to prove controllability of system  \eqref{eq:sistCont}-\eqref{eq: nao BG}.

A family $\F \subset \Ve V$ is called  {bracket-generating} if
\begin{equation*}
%\label{familiaBG}
  \Lie_p(\F)=T_p\,V, \qquad \forall p \in V
\end{equation*}
and a vector field  $X\in \Ve V$ is called compatible with a family
 $\F\subset \Ve V$ if
\begin{equation*}
    \cA_p(\F \cup \{X\})\subset \overline{\cA_p(\F)}, \quad p\in V.
\end{equation*}

For the proof of the following results see \cite{AgrSach01}.

\begin{proposition}
\label{thm: prop compativel}
Let $\F\subset \Ve V$ be a bracket-generating family. If a vector field $X\in \F$ is Poisson stable, then the vector field   $-X$ is compatible with $\F$.
\end{proposition}

\medskip

\begin{proposition}
\label{thm: prop controlavel}%\cite[corolario 8.6]{AgrSach01}
If $\F \subset \Ve V$ is a  {bracket-generating} family such that
the positive convex cone generated by $\F${ \begin{equation*}
        \mathrm{cone}(\F)=\left\{\sum_{i=1}^k a_i X_i: k \in \N, \, X_1,\dots, X_k\in \F, a_1, \dots, a_k\in C^\infty(V), a_i\geq0
    \right\}\subset \Ve V
\end{equation*}}
is symmetric, then $\F$ is controllable, i.e. $\cA_p(\F)=V$.
\end{proposition}

\medskip
The set $\mathrm{cone}(\F)$ being symmetric means that if $Y\in \mathrm{cone}(\F)$ then $-Y\in\mathrm{cone}(\F)$.

%%%%%%%%%%%%%%%%%%%%%%%%%%%%%
\section{Controllability analysis}
%%%%%%%%%%%%%%%%%%%%%%%%%%%%%
By means of Legendre transform, we convert Lagrangian system \eqref{eq: sis Lagrangiano} into  the Hamiltonian system
\begin{align*}
    &&\dot{p}_i=-\frac{\p H}{\p\theta_i}-m_iu(t), &&
    \dot{\theta}_i=\frac{\p H}{\p p_i}, \qquad i\in \{1,2\},&&
\end{align*}
where $H(\theta,p,t)=p\,\dot\theta-L(\theta,\dot\theta,t)$; this is a fourth dimensional (symmetric) system. Hamiltonian system state space is then $[0,\pi[\,\times[0,\pi[\,\times\R^2$.

For uncontrolled system, i.e. when  $u(t)\equiv 0$ for all $t>0$, Hamiltonian function is constant since it does not depend explicitly on time. On state space, solutions $(\theta,p)$ to Hamiltonian system are such that $H(\theta,p)\equiv c$, for different values of $c\in \R$ depending on initial conditions. These solutions remain on the compact defined by $H\leq c+1$. From Liouville theorem, the Hamiltonian system phase flow preserves volume and therefore, according to Theorem \ref{thm: Poicare}, the vector field associated with this flow is Poisson stable. Since the map $\Phi: (\theta, \dot \theta)\longrightarrow (\theta, p)$ is a diffeomorphism  and trajectories $(\theta, p)$ are Poisson stable on phase space, trajectories $(\theta, \dot \theta)$ on state space are Poisson stable.

Let $N$ be the set $H(\theta,p)\equiv c$, for different values of $c\in \R$, which is invariant under Hamiltonian dynamic and $M=\Phi^{-1}(N)$; then $M$ is invariant under Lagrangian dynamic.

\medskip
%%%%%%%%%%%%%%%%%%%%%%%%%%%
Two vector fields arise naturally associated with system \eqref{eq:sistCont},
\begin{equation}
\label{eq: campos originais}
f(z)\p_z \quad \text{and} \quad h(z)u(t)\p_z.
\end{equation}
{
On state space, the dynamic of uncontrolled system is governed by vector field $f\p_z$. We saw that trajectories $(\theta, \dot \theta)$ on state space are Poisson stable, therefore vector field $f\p_z$ is Poisson stable.}

Let
\begin{equation}
\label{eq: familia_F}
    \F=\{X_1,X_2,X_3,X_4\}
\end{equation}
where $X_1,X_2,X_3,X_4$ are vector fields obtained from vector fields in \eqref{eq: campos originais} and their Lie brackets
\begin{equation}
\label{eq: campos BG}
\begin{split}
&X_1=\Delta(\theta)\,f(z)\p_z= \Omega^T(z){\dth}+
a^T(z){\dom},\qquad
X_2=\Delta(\theta)\,h(z)\p_z=b^T(\theta){\dom},\\
&X_3=[X_1,X_2]=\bar{\Omega}^T(\theta){\dth}+
\bar{a}^T(\theta)\frac{\partial b}{\partial\theta}{\dom},\qquad
X_4=[X_2,X_3]=\bar{b}^T(\theta)\frac{\partial
b}{\partial\theta}{\dom}.
\end{split}
\end{equation}
Here the notation follows \eqref{eq: notacao 1}-\eqref{eq: notacao 4}.

\begin{theorem}
\label{thm: conj atin=orbita}
If $\F$ defined in \eqref{eq: familia_F} -\eqref{eq: campos BG} is a {bracket-generating} family of vector fields, then  $\cA_p(\F)=\cO_p(\F)$, for all $p\in
M$.
\end{theorem}

\medskip

\begin{pf}
Using vector fields in $\F$ we can obtain vector fields
$f(z)\p_z$ and $h(z)u(t)\p_z$. {Since $u(\cdot)$ is a scalar function both $h(z)|u(t)|\p_z$ and  $-h(z)|u(t)|\p_z$ belong to $\F$. Moreover, since $f(z)\p_z$ is Poisson stable, Proposition
\ref{thm: prop compativel} ensures that $-f(z)\p_z$ is compatible with $\F$. Therefore, for all $p\in M$,
$\cA_p(\F) =\cO_p(\F)$.\bbox}
\end{pf}

%%%%%%%%%%%%%%%%%%%%%%%%%%%%%%%%
\begin{corollary}
 Under conditions stated in Theorem {\ref{thm: conj atin=orbita}}, system \eqref{eq:sistCont}-\eqref{eq: nao BG}
is controllable.
\end{corollary}

\begin{pf}
If Theorem \ref{thm: conj
atin=orbita} holds, the cone
\begin{equation*}
    \mathrm{cone}(\F)=\left\{\sum_{i=1}^k a_i X_i: k \in \N, \, X_1,\dots, X_k\in \F, a_1, \dots, a_k\in C^\infty(V), a_i\geq0,
    \right\}
\end{equation*}
is symmetric and Proposition \ref{thm: prop controlavel} establishes that the system is controllable.\bbox
\end{pf}

Previous results assumed family $\F$ defined in \eqref{eq: familia_F}-\eqref{eq: campos
BG}
to be bracket generating. The proof that family $\F$ is bracket generating implies studying several subsets of manifold $M$ defined by means of analytic equations.
Although the proof is not hard it is  laborious. Main argument in the proof is that $\Lie_p(\F)$ has constant dimension
 for all  $p$ in some orbit of $\F$. Some details of the proof are omitted since only simple computations were involved.

%%%%%%%%%%%%%%%%%%%%%%%%%%%%%
\begin{lemma}
\label{thm: lema BG}
The family of vector fields $\F$ defined  in \eqref{eq: familia_F}-\eqref{eq: campos
BG} is   {bracket-generating}, i.e.
\begin{equation*}
    \Lie_p(\F)=T_p\,M, \qquad \forall p \in M
\end{equation*}
\end{lemma}

\begin{pf}
Let $\Gamma$ be the subset of $M$ on which vector fields $X_2$ and $X_4$ are linearly dependent and
$\Upsilon$  the subset of $M$ on which vector fields   $X_1$  and $X_3$ are linearly dependent on the direction $\dth$.
We denote by
\begin{equation*}
\overline{M}= M\backslash (\Gamma \cup \Upsilon)
\end{equation*}
the subset of $M$ on which vector fields defined in \eqref{eq: campos BG} are linearly independent. Then,
for all $ p\in \overline{M}$,
\begin{equation*}
    \dim \Lie_p(\F)=4.
\end{equation*}
Subsets $\Gamma$ and $\Upsilon$ are not $\F$ invariant, i.e. an orbit of $\F$ always leaves this sets.
In fact, although
$X_2$ and $X_4$ are linearly dependent over $\Gamma$ using $X_2$ and $X_3$ it is possible to leave  $\Gamma$ through an orbit of $\F$. Let  $x_2$ and $x_3$ denote the vector components for $X_2$ and $X_3$, respectively,
\begin{align*}
    x_2=\left(%
\begin{array}{c}
  0 \\
  0 \\
  b_1(\theta) \\
  b_2(\theta) \\
\end{array}%
\right),&&
x_3=\left(%
\begin{array}{c}
  \Delta(\theta)\,b_1(\theta) \\
  \Delta(\theta)\,b_2(\theta) \\
  *\\
   *\\
\end{array}%
\right).
\end{align*}
Since   $b_1(\theta)$ and $b_2(\theta)$ are not simultaneously zero, we get that $x_2$ and $x_3$ are not collinear.

Analogous situation occurs for $\Upsilon$, but now it is required to use also vector field $X_4$ in order to leave $\Upsilon$. To prove it, three main steps are needed. First we define a region $S_k$ and we solve equation $X_i\cdot S_k=0$. This last equation defines a set $S_{k+1}\subset S_k$ which is invariant  by $X_i$.\\

Let $\Sigma$ be the subset of  $M$ on which vector fields $X_1, X_3$ are linearly dependent on the direction  $\dom$ and $X_2,X_4$ are linearly dependent  (on direction $\dth$); $\Sigma$ is the intersection $\Gamma$ and $\Upsilon$. Here  vector fields $X_2$ and $X_3$ are linearly independent. Therefore $\Sigma$ has dimension two.\\

Let $q_0$ be a point in the orbit  $\cO_p(\F)$ through a point $p\in M$.
We study four different situations.

%%%%%%
First assume that $q_0 \in \Sigma$.
%%%%%%
Then, as it was seen earlier, the orbit $\cO_p(\F)$ leaves $\Sigma$. When
$\cO_{q_0}(\F)$ leaves $\Sigma$ either it enters $\overline{M},\,\Upsilon \setminus \Sigma$ or
$\Gamma\setminus \Sigma$.

%%%%%%
Now, assume that  $q_0 \in \overline{M}$.
%%%%%%
In $\overline{M}$ all vector fields in $\F$ are linearly independent. An
orbit of the family $\F$ through $q_0\in \overline{M}$ is, then, generated by four vector fields
and $\cO_{q_0}(\F)$ has dimension four. Therefore, the dimension of the tangent space to the orbit at any $z\in\cO_{q_0}(\F)$, $T_z \cO_{q_0}(\F)$, has dimension four. Moreover, by Hermann-Nagano theorem, for all $z \in\cO_{q_0}(\F)$,
\begin{equation*}
T_z\cO_{q_0}(\F)=\Lie_z (\F), \quad \text{and}\quad \dim \Lie_z \F=4.
\end{equation*}
Since for $z\in\overline{M}$,  the tangent space $T_z\overline{M}$ has dimension four we conclude that
\begin{equation*}
  T_z\overline{M}=\Lie_z(\F), \qquad \text{for}\,\, z\,\in  \overline{M}.
\end{equation*}

%%%%%%
Third situation occurs when $q_0 \in (\Upsilon \setminus \Sigma)$.
%%%%%%
The orbit $\cO_p(\F)$ leaves  $\Upsilon \setminus \Sigma$.
%%%%
Recall that in $\overline{M}$ the Lie algebra generated by $\F$ and evaluated at any $z\in\overline{M}$, $\Lie_z \F$, has dimension four and such dimension is constant for $z\in \overline{M}$.  If $\cO_{q_0}(\F)$ leaves  $\Upsilon$ to $\overline{M}$,  then for $z\in\Upsilon \setminus\Sigma$ the tangent space to the orbit $\cO_p(\F)$ must have dimension four:  $\dim T_z \cO_{q_0}(\F)=4$. And, once again,  from
Hermann-Nagano theorem
\begin{equation*}
\dim \Lie_z (\F)=4,\quad z\in (\Upsilon \setminus \Sigma).
\end{equation*}
This means that there must be four Lie brackets of vector fields in $\F$ that span $\Lie_z{\F}$,  $z\in (\Upsilon \setminus \Sigma)$.  On the other hand, the orbit  $\cO_{q_0}(\F)$ only leaves $\Upsilon \setminus
\Sigma$ to $\Gamma$ on the intersection of the two regions, that is, on $\Sigma$. But this is the first case studied.\\

%%%%%%
At last, let us assume that   $q_0 \in (\Gamma\setminus\Sigma)$.
%%%%%%
The orbit $\cO_{q_0}(\F)$ leaves $\Gamma\setminus\Sigma$.  If $\cO_{q_0}(\F)$
visits $\overline{M}$, by analogy with the second situation studied, there exist some Lie brackets spanning $\Lie_z  (\F)$, $z \in (\Gamma\setminus\Sigma)$.  On the other hand, the case when   $\cO_{q_0}(\F)$ visits
$\Upsilon \setminus \Sigma$ is the third case.\\

%%%%%%
Finally, we conclude that for $z\in M$, $T_z M=\Lie_z \F$ , i.e.  family $\F$ is {bracket-generating} in $M$.\bbox
\end{pf}
%%%%%%%%%%%%%%%%%%%%%%%%%%%%%
\section{Concluding remarks}
%%%%%%%%%%%%%%%%%%%%%%%%%%%%%%%%
{
In this paper we considered a mechanical system consisting in two linked planar inverted pendulums on a cart. The main result obtained is the proof of existence of a bracket-generating family of vector fields $\F$ for which the Lie algebra has constant dimension over each orbit of the family through a point. In the proof of Lemma \ref{thm: lema BG} it is necessary to guarantee \eqref{eq: nao BG}. If equality holds in  \eqref{eq: nao BG} it is not  possible to prove that the orbit of family $\F$ leaves $\Sigma$ and therefore to prove that $\F$ is bracket-generating.\\
}

As vector fields $X_2$ and $X_4$ in \eqref{eq: campos BG} only have non zero components in two directions,
$\dom$, it seams natural to look for two other vector fields with non zero components only in directions  $\dth$. For the subset of $M$ where $X_2,X_4$ are linearly independent these vector fields form a basis for directions $\partial \omega$.  Therefore, it is possible to span new vector fields, say $Y_1$ and $Y_3$, which components associated to directions $\partial \omega$ are zero
\begin{align}
\begin{split}\label{eq:def_Y}
Y_1&= X_1-(\gamma_2(z)\,X_2+\gamma_4(z)\, X_4), \qquad \gamma_2,\,\gamma_4\in\R\\
    Y_3&= X_3-(\tilde\gamma_2(z)\,X_2+\tilde\gamma_4(z)\, X_4), \qquad
    \tilde\gamma_2,\,\tilde\gamma_4 \in\R.
\end{split}
\end{align}
That is, new vector fields are such that   $Y_1=\Omega^T(z){\dth}$ and
$Y_3=\bar{\Omega}^T(\theta){\dth}$  given that (\ref{eq:def_Y}) does not alters components on directions  $\dth$. This way, we could had chosen to work with another family of vector fields, say
\begin{equation*}
    \mathcal{G}=\{Y_1,X_2,Y_3,X_4\}.
\end{equation*}
Since vector fields in $\mathcal{G}$ result from linear combinations of vector field in $\F$, this new family is also bracket generating family for system  \eqref{eq:sistCont}. Yet, vector fields in  $\mathcal{G}$ are not defined for all points in $M$ so it would be necessary to use other vector fields to study $\Gamma$.\\

%%%%%%%%%%%%%%%%%%%%%%%%%%%%%
\section{Acknowledgements}
%%%%%%%%%%%%%%%%%%%%%%%%%%%%%%%%
The author thanks A.V. Sarychev for first introducing her to the problem and for several enlightening
discussions.
%%%%%%%%%%%%%%%%%%%%%%%%%%%%%%%%

%%%%%%%%%%%%%%%%%%%%%%%%%%%%%%%%%%%%%%
\end{document}